\documentclass[a4paper,english,leqno, 10pt]{amsart}
\usepackage{amssymb,amsfonts,amsmath,amsthm}
\usepackage{epsfig,epsf,curves,epic,eepic}
\usepackage{enumerate}
\usepackage{bbm}
\usepackage[mathscr]{eucal} 
\newtheorem{thm}{Theorem}
\newtheorem{lem}{Lemma}

\newtheorem{cor}{Corollary}
\theoremstyle{remark}
\newtheorem{exmp}{Example}
\newtheorem{defn}{Definition}

\newtheorem{rem}{Remark}
\newcounter{fig}
\numberwithin{equation}{section}

\newcommand{\OR}{\overline{\mathbb{R}}}
\newcommand{\N}{\mathbb{N}}
\newcommand{\U}{\mathbb{U}}
\newcommand{\C}{\mathbb{C}}
\newcommand{\OC}{\overline{\mathbb{C}}}

\newcommand{\cc}{\mathcal{C}}
\newcommand{\cd}{\mathcal{D}}

\newcommand{\cf}{\mathcal{F}}
\newcommand{\cj}{\mathcal{J}}
\newcommand{\cs}{\mathcal{S}}
\DeclareMathOperator{\Aut}{Aut}

\DeclareMathOperator{\diam}{diam}
\DeclareMathOperator{\poles}{poles}
\DeclareMathOperator{\zeroes}{zeroes}
\DeclareMathOperator{\lcm}{lcm}
\DeclareMathOperator{\spec}{spec}
\DeclareMathOperator{\cells}{cells}

\newcommand{\V}[1]{\mathit{V#1}}
\newcommand{\E}[1]{\mathit{E#1}}
\begin{document}
$\mbox{}$\vskip-0.5cm
\title[Green functions on self-similar graphs and bounds for the spectrum]{Green functions on self-similar graphs and\\ bounds for the spectrum of the Laplacian}
\author[B.\ Kr\"on]{B.\ Kr\"on$^{\textstyle{\,\star}}$}
\begin{abstract}
Combining the study of the simple random walk on graphs, generating functions (especially Green functions), complex dynamics and general complex analysis we introduce a new method for spectral analysis on self-similar graphs.\par
First, for a rather general, axiomatically defined class of self-similar graphs a graph theoretic analogue to the Banach fixed point theorem is proved. The subsequent results hold for a subclass consisting of `symmetrically' self-similar graphs which however is still more general then other axiomatically defined classes of self-similar graphs studied in this context before: we obtain functional equations and a decomposition algorithm for the Green functions of the simple random walk Markov transition operator $P$. Their analytic continuations are given by rapidly converging expressions. We study the dynamics of a probability generating function $d$ associated with a random walk on a certain finite subgraph (`cell-graph'). The reciprocal spectrum $\spec^{-1}\!P=\{1/\lambda\mid \lambda\in\spec P\}$ coincides with the set of points $z$ in $\OR\setminus (-1,1)$ such that there is Green function which cannot be continued analytically from both half spheres in $\OC\setminus\OR$ to $z$. The Julia set $\cj$ of $d$ is an interval or a Cantor set. In the latter case $\spec^{-1}\!P$ is the set of singularities of all Green functions. Finally, we get explicit inner and outer bounds, $\cj\subset\spec^{-1}\!P\subset\cj\cup\cd,$ where $\cd$ is the set of the $d$-backward iterates of a finite set of real numbers.
\end{abstract}

\bibliographystyle{plain}
\thanks{$^{\textstyle{\star}}$ The author was supported by the projects Y96-MAT and P14379-MAT of the Austrian Science Fund. Current address: Erwin Schr\"odinger Institute (ESI) Vienna, Boltzmanngasse 9, 1090 - Wien, e-mail: bernhard.kroen@univie.ac.at.}
\thanks{This work is part of the author's PhD-thesis \cite{kroen01spectral}.}
\thanks{Mathematics Subject Classification 60J10, (30D05, 05C50)}
\maketitle

\section{Introduction}
Self-similar graphs can be seen as discrete versions of fractals (more precisely: compact, complete metric spaces defined as the fixed set of an iterated system of contractions, see Hutchinson in \cite{hutchinson81fractals}). The simple random walk is a crucial tool in order to study diffusion on fractals, see Barlow and Kigami \cite{barlow97localized}, Barlow  and Perkins \cite{barlow88brownian}, Grabner \cite{grabner97functional2}, Lindstr{\o}m \cite{lindstroem90brownian} and many others. In the present paper we study Green functions and the spectrum of the corresponding Markov transition operator $P$ and (equivalently) the Laplacian $\Delta=I-P$, where $I$ is the identity.\par
There are different methods of spectral analysis on self-similar structures. Barlow and Kigami \cite{barlow97localized} and Sabot \cite{sabot00pure} used localized eigenfunctions to study the spectrum and to prove that it is a pure point spectrum. A first heuristic result on explicit spectra is due to Rammal \cite{rammal84spectrum} who investigated the spectrum of the Sierpi\'{n}ski graph in the setting of statistical physics. Malozemov and Teplyaev studied the spectral self-similarity of operators in a series of papers. Their functional analytic method is fundamentally different from our probabilistic approach. Malozemov computed the spectrum of the Koch graph in \cite{malozemov93integrated} and the spectrum of a self-similar tree in \cite{malozemov95random}. The spectrum of the Sierpi\'{n}ski graph was computed by Teplyaev in \cite{teplyaev98spectral}. In \cite{malozemov95pure} they discussed spectral properties of a class of self-similar graphs which can be seen as discrete analogues to finitely ramified fractals such that the boundaries of their cells consist of two points. We give a first axiomatic definition of self-similar graphs corresponding to fractals with an arbitrary number of boundary points. With our method we do not only obtain information about the spectrum but we are also able to describe the analytic properties of the Green functions precisely. This axiomatic definition of self-similarity of graphs and the central results of the present paper are part of the author's PhD thesis \cite{kroen01spectral}
\footnote{Results and construction of self-similar graphs were presented at the DMV-conference (Dresden, September 2000), `Dynamic Odyssey'  (Luminy-Marseille, February 2001), `Fractals' (Graz, June 2001) and `Random walks and Geometry' (Vienna, June/July 2001).}.
In a recent pre-print \cite{malozemov01self} Malozemov and Teplyaev constructed a similar but more restricted class of self-similar graphs and they obtained the same bounds for the spectrum as in Theorem \ref{thm:bounds} of this paper.\par
Bartholdi and Grigorchuk, see \cite{bartholdi00croissance} and \cite{bartholdi00spectrum}, have computed the spectra of several Schreier graphs of fractal groups of intermediate growth. The reader is also referred to the survey article of Bartholdi, Grigorchuk and Nekrashevych \cite{bartholdi02fractal}.\par
We introduce a new class of rather general self-similar graphs in Section \ref{ssg}. In this definition there occurs a function $\phi$ which maps a set of vertices to the set of all vertices. This function is a contraction with respect to the natural graph metric and it can be interpreted as the `self-similarity function' of the graph. A graph theoretic analogue of the Banach fixed point theorem is proved: Either $\phi$ fixes exactly one origin vertex, or it contracts towards exactly one `cell' of the graph. The simple random walk on self-similar graphs is recurrent if the graph has bounded geometry, which means that the set of vertex degrees is bounded, the cells are finite and the numbers of vertices in the boundaries of the cells are bounded. Properties concerning bounded geometry and volume growth of self-similar graphs were studied by the author in \cite{kroen02growth}.\par
\begin{figure}
\refstepcounter{fig}\label{vicsek}
\includegraphics{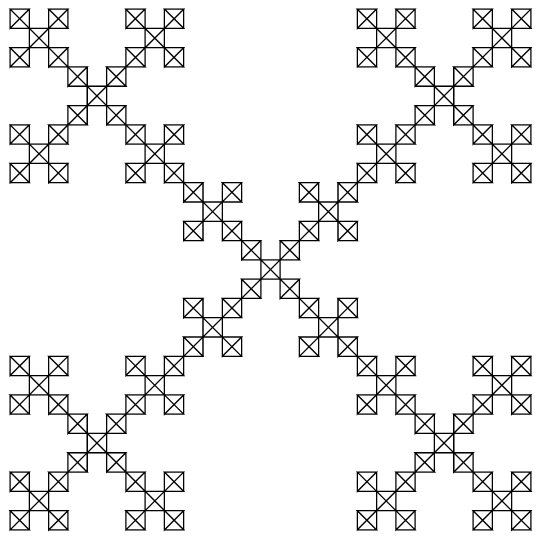}
\begin{center}
\emph{Figure \ref{vicsek}}
\end{center}
\end{figure}

In Section \ref{symmetric} the class of `symmetrically self-similar' graphs is defined. They correspond to finitely ramified fractals for which the renormalization problem for the simple random walk (see Lindstr{\o}m \cite{lindstroem90brownian} or Kigami \cite{kigami93harmonic}) can be solved by using only one variable. In other words: Starting at a point in the boundary of a cell in a fractal, the transition times for the Brownian motion to hit any other point in this boundary have the same distribution for all boundary points. Our class of graphs contains the graphs studied by Malozemov and Teplyaev but also other well known graphs such as the extensively studied Sierpi\'{n}ski graph, see Barlow and Perkins \cite{barlow88brownian}, Grabner and Woess \cite{grabner97functional1}, Hambly \cite{hambly00asymptotics} and Jones \cite{jones96transition} or the Vi\v{c}ek graph which corresponds to the Vi\v{c}ek snowflake, see Metz and Hambly \cite{hambly98homogenization} and \cite{metz93how}. A finite part of this graph can be seen in Figure \ref{vicsek}.

\par
We apply a variable substitution technique of generating functions (see Goulden and Jackson \cite{goulden83combinatorial}) in Section \ref{ssg+gf}. This technique was used by Grabner and Woess in \cite{grabner97functional1} to derive a functional equation for a certain Green function on the Sierpi\'{n}ski graph. Rammal and Toulouse used physical arguments to obtain such a functional equation, see \cite{rammal84random} and \cite{rammal83random}. Grabner and Woess calculated the asymptotic fluctuation behaviour of the $n$-step return probabilities. Grabner gave further applications of this substitution method in the context of the analysis of stopping times of Brownian motion on the Sierpi\'{n}ski gasket, see \cite{grabner97functional2}. We state similar functional equations for all Green functions on any symmetrically self-similar graph. Such functional equations reflecting the self-similarity of a graph, a group or a fractal occur often in the literature. We refer for example to the papers mentioned above: Bartholdi, Grigorchuk and Nekrashevych (\cite{bartholdi00croissance}, \cite{bartholdi00spectrum} and \cite{bartholdi02fractal}) and Malozemov and Teplyaev (\cite{malozemov93integrated}, \cite{malozemov95random}, \cite{malozemov95pure} and  \cite{teplyaev98spectral}).\par
In Section \ref{decompose} the functional equations for the Green functions of the whole graph are extended to a decomposition formula. This formula together with the fixed point theorem of Section \ref{ssg} implies that any Green function can be written as a finite sum of terms consisting of products and concatenations of finitely many Green functions (of the origin cell and the origin vertex) together with the \emph{return function} $f$, the \emph{transition function} $d$ and the \emph{inner transition functions} $h$. The functions $f$, $d$ and $h$ are rational Green functions of the simple random walk with absorbing boundary on the finite subgraph which is spanned by a cell. This decomposition is one of the central tools in studying the singularities of the Green functions (see also \cite{kroen02asymptotics}) and in computing the spectrum of the Laplacian.\par
Green functions of the simple random walk with absorbing boundaries on finite subgraphs spanned by the so-called `$n$-cells' are studied in Section \ref{ncells}. These rational functions can be decomposed into the return function, the transition function and the inner transition functions. The zeroes of the return function turn out to be poles of the inner transition functions.\par
Preliminary results concerning the dynamics of the transition function $d$, the analytic continuations of the Green functions and the spectrum of the transition operator are stated in Section \ref{dynamics}. We define $\cd$ as the set of $d$-backward iterates of the poles of the return function and the inner transition functions. This is the set of poles of the Green functions studied in Section \ref{ncells} being contained in $\OR\backslash (-1,1)$. The reduction algorithm in Section \ref{decompose} yields the analytic continuation of all Green functions to the immediate basin of attraction $A_0^*$ of the fixed point $z=0$ of $d$ except for some points in $\cd$. It will turn out that $A_0^*$ coincides with the Fatou set of $d$.\par
In Section \ref{notisolated} we put all these results together and give explicit bounds for the spectrum:
The reciprocal spectrum
\[\spec^{-1}\!P=\{1/\lambda\mid \lambda\in\spec P\}\]
of the transition operator $P$ of the simple random walk is the set of real numbers $z$ in $\OR$ such that the Green functions cannot be continued uniquely analytically to $z$ from both half spheres in $\OC\backslash\OR$ and
\[\cj\ \subset\ \spec^{-1}\!P\ \subset\ \cj\cup\cd\ \subset\ \OR\backslash (-1,1).\]
A crucial role is played by the singularity $z=1$ of the Green functions.\par
The last Section \ref{julia} is devoted to the Julia set $\cj$. The Julia set is either a real Cantor set or it is a real interval in the extended complex plane $\OC$. If $\cj$ is a Cantor set, then the reciprocal spectrum coincides with the set of all singularities of all Green functions. We give the example of the two sided infinite line, where the reciprocal spectrum is $\OR\backslash (-1,1)$. In this case the Green function has a branch cut along $\OR\setminus(-1,1)$. Here the transition functions are conjugated to the Chebychev polynomials.

\section{Self-similar graphs}\label{ssg}

Graphs $X=(\V X,\E X)$ are always assumed to be connected, locally finite, infinite and without loops or multiple edges. If $X$ is undirected then each $e=\{x,y\}$ is a set of vertices $x$ and $y$ in $\V X$, if $X$ is directed, then $\E X\subset\V X\times\V X$. A \emph{path of length $n$} from $x$ to $y$ is an $(n+1)$-tuple of vertices
\[(x_0=x,x_1,\ldots,x_n=y)\]
such that $x_i$ is adjacent to $x_{i+1}$ for $0\le i\le n-1$. The distance $d_{X}(x,y)$ is the length of the shortest path from $x$ to $y$. We define the \emph{vertex boundary} or \emph{boundary} $\theta C$ of a set $C$ of vertices in $\V X$ as the set of vertices in $\V X\backslash C$ which are adjacent to a vertex in $C$. The \emph{closure} of $C$ is defined as $\overline{C}=C\cup \theta C$. Since we will use the topological closure only for sets of complex numbers, this notion will cause no misunderstandings. The set $C$ is called \emph{connected} if every pair of vertices in $C$ can be connected by a path in $X$ that does not leave $C$. The set $C$ is \emph{adjacent} to a vertex $x$ if $x\in\theta C$. Two sets of vertices are \emph{adjacent} if their boundaries are not disjoint.\par
For a set $F\subset \V X$ let $\cc_{X}(F)$ denote the set of connected components in $\V X\backslash F$. We define a new graph $X_{F}$ by setting $\V X_{F}=F$ and connecting two elements $x$ and $y$ in $\V X_{F}$ by an edge if and only if there exists a $C\in\cc_{X}(F)$ such that $x$ and $y$ are in the boundary of $C$. The graph $X_{F}$ is then called the \emph{reduced graph of $X$} with respect to $F$.

\begin{defn}\label{def:ssg}
A graph $X$ is called \emph{self-similar} with respect to $F$ if
\begin{enumerate}[(F1)]
\item no vertices in $F$ are adjacent in $X$,
\item the intersection of the closures of two different components in $\cc_{X}(F)$ contains no more than one vertex and
\item there exists an isomorphism $\psi: X\to X_{F}$.
\end{enumerate}

We will also write $\phi$ instead of $\psi^{-1}$ and $F^{n}$ instead of $\psi^{n}F$, $F^0$ is defined as $\V X$. The components of $\cc_{X}(F^{n})$ are called \emph{$n$-cells}, 1-cells are also just called \emph{cells}. Furthermore we define functions
\[S: \cc_{X}(F)\to \cc_{X}(F^{2})\]
where $S(C)$ is the 2-cell that contains the cell $C$ and
\[\phi_{S}: \cc_{X}(F)\to \cc_{X}(F)\]
such that $\phi_{S}C$ is the cell whose boundary is $\phi\,\theta S(C)$; equivalently: $\phi_{S}C=\phi (S(C)\cap F)$. A cell $C$ such that $\phi_{S}C=C$ is called an \emph{origin cell}, a vertex which is a fixed point of $\phi$ is called an \emph{origin vertex}. If $\phi_{S}^{m}C=C$ for some integer $m$ then we call $C$ a \emph{periodic origin cell with period $m$}.
\end{defn}

Note that the boundary of any cell in $X$ spans a complete graph in $X_F$.

\begin{lem}\label{twocells}
Let $C$ be an $m$-cell and $D$ be an $n$-cell. Then either $m=n$ and $C=D$ or
\[|\theta C \cap \theta D|\le 1.\]
\end{lem}

\begin{proof}
If $m=n$, then this is equivalent to Axiom (F2) for the graph $X_{F^{m-1}}=X_{F^{n-1}}$. Otherwise let us suppose that $m<n$. If there are two different elements $x$ and $y$ in $\theta C \cap \theta D$, then they are adjacent in $X_{F^m}$. The vertices $x$ and $y$ are elements of $F^n$ which is a subset of $F^{m+1}$. Since the graph $X_{F^m}$ is self-similar with respect to $F^{m+1}$, this contradicts Axiom (F2) for $X_{F^m}$.
\end{proof}

The following theorem can be seen as a graph theoretic analogue for the Banach fixed point theorem.

\begin{thm}\label{banach}
Let $X$ be a self-similar graph with respect to $F$. Then $X$ is also self-similar with respect to $F^{n}$. Either there exists exactly one origin cell $O$ and for every cell $C$ there exists an integer $n$ such that
\[\phi_{S}^{n}C=O\]
or $\phi$ has exactly one origin vertex and for any cell $C$ there exists a periodic origin cell $O_{m}$ with period $m$ and an integer $n$ such that
\[\phi_{S}^{n}C=O_{m}.\]
\end{thm}

\begin{proof}
The graphs $X$ and $X_{F}$ are isomorphic, and $X_{F}$ is isomorphic to $X_{F}$ reduced by $\psi(F)$ which implies that $X$ is isomorphic to $X_{F^2}$. The first statement now follows by induction.\par
Let $x$ and $y$ be vertices in $F$ with $d_{X_{F}}(x,y)=n$. Since no vertices in $F$ are adjacent in $X$, every path in $X$ connecting these vertices must have length at least $2n$. Thus 
\[d_{X}(x,y)\ge 2\cdot d_{X}\big(\phi x,\phi y\big)=2\cdot d_{X_{F}} (x,y)\]
for all pairs of vertices $x$, $y$ in $F$. For any cell $C$ in $\cc_{X}(F)$ 
\[d_{X}(\theta C,\theta\phi_{S}C)\ge d_{X}\big(\theta S(C),\theta S(\phi_{S}C)\big)\ge 2\cdot d_{X}\big(\phi\,\theta S(C),\phi\,\theta S(\phi_{S}C)\big)= 2\cdot d_{X}\left(\theta \phi_{S}C,\theta \phi_{S}^{2}C\right)\]
and the sequence $\left(d_{X}(\theta \phi_{S}^{n-1}C,\theta \phi_{S}^{n}C)\right)_{n\in \mathbb{N}}$ must eventually be zero.\par
If $\phi$ has fixed points $x$ and $y$, then
\[d_{X}(x,y)\ge 2\cdot d_{X}\big(\phi x,\phi y\big)= 2\cdot d_{X}(x,y)\]
and $x$ equals $y$.\par
Let $C_{1}$ and $C_{2}$ be different cells such that $\phi_{S}C_{1}=C_{2}$ and $d_{X}(\theta C_{1},\theta C_{2})=0$. By Axiom (F2) there is exactly one vertex $x$ in $\theta C_{1}\cap\theta C_{2}$. In the case that $C_{1}$ and $C_{2}$ are in the same 2-cell, $C_{2}$ is an origin-cell. Otherwise we distinguish: If $\phi\,\theta C_{2}$ lies in $S(C_{2})$, then $\phi (S(C_{2})\cap F)$ is the cell $\phi_S C_2$ which contains $\phi\,\theta C_{2}$, and $\phi_S C_2$ is an origin cell. If $\phi\,\theta C_{2}$ does not lie in $S(C_{2})$, then $\theta S(C_{2})$ and $\theta \phi_{S}C_{2}$ must have non-empty intersection, because $d_{X}(\theta C_{2},\theta \phi_{S}C_{2})=0$. By Lemma \ref{twocells} this intersection consists only of the vertex $x$. Thus
\[x\in \theta C_{1}\cap \theta \phi_{S}C_{1}\cap\theta \phi^{2}_{S}C_{1},\]
and proceeding by induction we obtain $x\in\theta\phi^{n}_{S}C_{1}$ for every integer $n$. Since $X$ is locally finite there must be indices $k$ and $m$ such that $\phi_{S}^{k+m}C_{1}=\phi_{S}^{k}C_{1}$. The closures of all different 2-cells containing one of the cells with $x$ in their boundaries must have only the vertex $x$ in common. The function $\phi$ maps these 2-cells onto the cells having $x$ in their boundaries and therefore it must fix $x$.\par
Suppose there is an origin cell $C$ and no fixed point of $\phi$. If $\theta C$ and $\theta S(C)$ are disjoint, then no cell $D$ which is not in $S(C)$ can be fixed by $\phi_S$, because $d_X(\theta \phi^n_S D,\theta \phi^n_S C)$ must be zero for sufficiently large integers $n$, and $C$ is the only origin cell. If there is a vertex $x$ in $\theta C\cap \theta S(C)$, then by the above arguments
\[x\in \phi^{n}_{S}C\]
for any integer $n$ and $x$ must be fixed by $\phi$.
\end{proof}

\begin{lem}\label{components}
Let $X$ be self-similar with respect to $F$ and the corresponding isomorphism $\phi$ and let $o$ be an origin vertex. Then there is an integer $k$ such that the subgraphs $\hat A$ which are spanned by $\overline A=A\cup \{o\}$ for components $A$ in $\cc_X(\{ o \} )$ are self-similar graphs with respect to
\[F^k\cap \overline A \text{\quad and\quad} \phi^k|_{F^k\cap \overline A},\]
and they have exactly one origin cell $O_A$.
\end{lem}

\begin{proof}
The fixed point $o$ of $\phi$ must be contained in the boundary of every periodic origin cell. Since $X$ is locally finite, there exist only finitely many periodic cells $C_1$, $C_2,\ldots,C_r$ with periods $m_{1}$, $m_{2},\ldots,m_{r}$. For
\[k=\lcm (m_{1}, m_{2},\ldots,m_{r})\]
we have $\phi_S^{k}C_{m_{i}}=C_{m_{i}}$, where $i\in\{1,2,\ldots,r\}$, and all cells $C_i$ are origin cells with respect to $F^k$ and $\phi^k$. Let $A_{i,n}$ be the closed $n$-cell which contains $C_i$. Then, for $n\ge 2$,
\[\phi^k (A_{i,n}\cap F^k)=A_{i,n-1}\]
and if $C_i$ is the origin cell with respect to $F^k$ which lies in the component $A$ in $\cc_X(\{ o \} )$, then
\[\bigcup_{n=1}^{\infty}A_{i,n}=\overline A.\]
\end{proof}

\begin{lem}\label{recurrence}
Let all cells, $\sup \{|\theta C|\mid C\in\cc_X(F)\}$ and $\sup \{\deg_X x\mid x\in F\}$ be finite. Then the simple random walk on $X$ is null recurrent.
\end{lem}

\begin{proof}
We apply Nash-Williams' recurrence criterion, see \cite{nash59random}, as a special case of Rayleigh's short-cut method which, for example, can be found in the book of Doyle and Snell \cite[Chapter~6]{doyle84random}, or in the book of Woess \cite[Theorem~2.19~and~Corollary~2.20]{woess00random}.\par
Let $o$ be an origin vertex. By Lemma \ref{components} we can choose $F$ such that all cells $C$ with $o\in\theta C$ are origin cells. Let $O_n$ be the union of the closures of all $n$-cells $C_n$ with $o\in\theta C_n$. If $X$ has no origin vertex, we define $O_n$ as the closure of the $n$-cell which contains the (in this case unique) origin cell, see Theorem \ref{banach}. The sequence
\[(A_n)_{n\in\N},\quad A_1=O_1 \mbox{\quad and\quad} A_{n+1}=O_{n+1}\backslash O_n,\]
is a one-dimensional partition of $\V X$ in the following sense:
\[\bigcup_{n=1}^{\infty}A_n=\V X,\]
and for any $n\ge 2$ the sets $A_n$ are finite and the vertices in $A_n$ are only adjacent to vertices in $A_{n-1}$, $A_n$ or $A_{n+1}$. Since
\[\sup \{|\theta C|\mid C\in \cc_X(F)\}=\sup \{|\theta C|\mid C\in \cc_X(F^n)\}\]
for every $n\in\N$ and $\{\deg_X (x)\mid x\in F\}$ is bounded, the number $a_n$ of edges connecting two sets $A_n$ and $A_{n+1}$ is also bounded. Thus
\[\sum_{n=1}^\infty\frac{1}{a_n}=\infty\]
and, by the Nash-Williams' recurrence criterion, the simple random walk on $X$ is recurrent. Recurrence of the simple random walk on infinite graphs is always null recurrence, see for example \cite[Theorem~1.18]{woess00random}.
\end{proof}

\section{Symmetrically self-similar graphs}\label{symmetric}

\begin{defn}
Let $X$ be a connected graph which is self-similar with respect to $F$ and et $\hat C$ be the subgraph of $X$ which is spanned by the closure $\overline C$ of a cell $C$. From now on we suppose that
\begin{enumerate}
\item[(S1)] all cells are finite and for any pair of cells $C$ and $D$ in $\cc_{X}(F)$ there exists a graph isomorphism $\alpha:\hat{C}\mapsto\hat{D}$ such that $\alpha\theta C=\theta D$.
\end{enumerate}\par
The graph $X$ is called \emph{simply symmetrically} self-similar if
\begin{enumerate}
\item[(S2)] $\Aut(\hat{C})$ acts transitively on $\theta C$,
\end{enumerate}\par
and \emph{doubly symmetrically} self-similar or just \emph{symmetrically} self-similar if
\begin{enumerate}
\item[(S3)] $\Aut(\hat{C})$ acts doubly transitively on $\theta C$, which means that it acts transitively on the set of ordered pairs
\[\left\{ (x,y)\mid x,y \in \theta C,\ x\ne y \right\},\]
where $g((x,y))$ is defined as $(g(x),g(y))$ for any $g\in \Aut(\hat{C})$.
\end{enumerate}
Let $\cells_X x$ be the number of cells which are adjacent to $x$. If $X$ satisfies (S1), we write $\theta_X$ for the number of vertices in the boundary $\theta C$ of some cell $C$. 
\end{defn}

If $X$ satisfies Axiom (S1), then the number $\mu$ of cells in a 2-cell is independent of the choice of this 2-cell and for any cell $C$ the graph $\hat C$ consists of $\mu$ copies of the complete graph $K_{\theta_X}$. This fact and many more details can be found in \cite{kroen02growth}.

\vskip4mm
\begin{center}
\begin{picture}(60,60)
\put(-20,25){$\hat C$}
\put(0,0){\circle*{4}}
\put(0,60){\circle*{4}}
\put(60,0){\circle*{4}}
\put(60,60){\circle*{4}}
\put(20,0){\circle*{2.5}}
\put(40,0){\circle*{2.5}}
\put(20,60){\circle*{2.5}}
\put(40,60){\circle*{2.5}}
\put(0,20){\circle*{2.5}}
\put(20,20){\circle*{2.5}}
\put(40,20){\circle*{2.5}}
\put(60,20){\circle*{2.5}}
\put(0,40){\circle*{2.5}}
\put(20,40){\circle*{2.5}}
\put(40,40){\circle*{2.5}}
\put(60,40){\circle*{2.5}}
\drawline(0,0)(60,60)(40,60)(40,0)(60,0)(0,60)(20,60)(20,0)(0,20)(60,20)(60,0)(40,0)(60,20)
\drawline(20,0)(0,0)(0,20)
\drawline(0,60)(0,40)(20,60)(0,40)(60,40)(60,60)(60,40)(40,60)
\refstepcounter{fig}\label{vicell}
\end{picture}
\vskip3mm
\emph{Figure \ref{vicell}}
\end{center}

Figure \ref{vicell} shows the subgraph $\hat C$ for the Vi\v{c}ek graph. It consists of five copies of the complete graph $K_4$. The vertices which are drawn fat belong to $F$; they constitute the boundary of $C$. Note that there are vertices which are only contained in the boundary of one cell, see Figure \ref{vicsek}, no matter whether the graph has an origin vertex or an origin cell. 

\begin{lem}\label{cells}
Let $X$ be a self-similar graph satisfying Axiom (S1). Then for any $v\in F$
\[\cells_{X}v\cdot(\theta_X-1)=\deg_{X_{F}} v.\]
\end{lem}

\begin{proof}
Each $w$ in the boundary of a cell $C$ with $v\in\theta C$ which is different from $v$ corresponds to a neighbour of $ v$ in $X_F$.
\end{proof}

\begin{lem}\label{degree}
Let $X$ be a simply symmetrically self-similar graph. Then the following conditions are equivalent:
\begin{enumerate}
\item $X$ has bounded geometry.
\item If $v$ is in the boundary of some cell $C$, then $\theta_X-1= |\{x\in C\mid x\sim v\}|$.
\item $\deg_X v=\deg_{X_F} v$ for all $v\in F$.
\end{enumerate}
\end{lem}

\begin{proof}
$(1)\,\Rightarrow\,(2)$. Suppose $X$ has bounded geometry but
\begin{equation}\label{a_eqn}
\frac{\deg_X v}{\cells_X v}=|\{x\in C\mid x\sim v\}|= a(\theta_X-1)
\end{equation}
for some $a\ne 1$. By the simple symmetry this number is independent of the choice of $v$ and $C$. Since $X$ and $X_F$ are isomorphic with respect to $\psi$, by \eqref{a_eqn} and Lemma \ref{cells}, we have
\[\deg_{X_F}\psi v=\deg_X v=\cells_X v\cdot a(\theta_X-1)=a\cdot \deg_{X_F}v\]
and for any positive integer $n$
\[\deg_{X_F} \psi^n v=a^n\cdot \deg_{X_F} v.\]
Thus either $\big(\deg_{X_F} \psi^n v\big)_{n\in\N}$ tends to zero, which is impossible, or $X_F$, and therefore also $X$, do not have bounded geometry.\par
$(2)\,\Rightarrow\,(3)$. If condition (2) is satisfied, then
\[\deg_X v=\cells_X v \cdot (\theta_X-1)=\deg_{X_F} v\]
for all $v\in F$.\par
$(3)\,\Rightarrow\,(1)$.
The condition
\[\deg_X v=\deg_{X_F}v\]
implies
\[\deg_X v=\deg_{X_{F^n}}v\]
for all vertices $v$ in $F^n$. By Axiom (F3) we have
\[\deg_X x=\deg_{X_F}\psi x=\deg_{X_{F^n}}\psi^n x\]
for all $x$ in $\V X$. For $x=\phi^n v$ this means
\[\deg_X\phi^n v=\deg_{X_{F^n}}\psi^n \phi^n v=\deg_{X_{F^n}} v=\deg_X v.\]
Because of Axiom (S1) the degrees of the vertices in $\V X\backslash F$ are bounded. For each $v\in F$ there is an $n$ such that $\phi^n v$ is either a vertex in $\V X\backslash F$ or an origin vertex $o$. Thus, for any cell $C$,
\[\max \big\{ \deg_X x \mid x\in \{o\}\cup C\big\}\]
is the maximal degree in $X$.
\end{proof}

\section{A functional equation for Green functions\\ on symmetrically self-similar graphs}\label{ssg+gf}

Let $P=(p(x,y))_{x,y\in \V X}$ be the matrix of the transition operator of the simple random walk on the space of functions
\[\ell^{2}(X)=\{f:\V X\to \mathbb{C}\mid\sum_{x\in \V X}\deg_X x\cdot |f(x)|^{2}<\infty\}\]
with inner product
\[(f,g)=\sum_{x\in \V X}\deg_X x\cdot f(x)\overline{g(x)}.\]
Throughout this article, we consider the compact Riemannian sphere $\OC$, the closure of the complex numbers $\C$. The Green function for vertices $x$ and $y$ in $\V X$ is defined as the generating function
for the $n$-step transition probabilities $p^{(n)}(x,y)$ from $x$ to $y$. Let $\U$ denote the complex, open unit disc, then
\[G(x,y|z)=\sum_{n=0}^{\infty}p^{(n)}(x,y)z^{n},\qquad z\in \U.\]
In matrix notation
\[\big(G(x,y|z)\big)_{x,y\in \V X}=\sum_{n=0}^{\infty} P^n z^n= (I-zP)^{-1},\]
whenever $1/z$ is not in the $\ell^2$-spectrum of $P$.\par
Following Grabner and Woess in \cite{grabner97functional1} we use the combinatorics of paths to derive a functional equation for the Green functions. Let $X$ be symmetrically self-similar. We write $\Pi_{X}(x,y)$ for the set of paths in $X$ from $x$ to $y$ and $\Pi_{X}^{*}(x,y)$ for those paths from $x$ to $y$ which meet $y$ only at the end. The \emph{weight} of a path $\pi=(x_0,\ldots,x_n)$ is defined as
\[W(\pi|z)=\prod_{i=0}^{n-1}\frac{z}{\deg_X x_{i}}\]
where $z\in\C$. For a set of paths $\Pi$ we set $W_{X}(\Pi|z)=\sum_{\pi\in\Pi}W_{X}(\pi|z)$. Then we have
\[G(x,y|z)=W_{X}(\Pi_{X}(x,y)|z).\]
For vertices $v$ and $w$ in $F$ and a path $\pi=(v\!=\!x_{0},x_{1},\ldots,w\!=\!x_{n})$ in $X$ we define $\tau_{j}=\tau_{j}(\pi)$ by $\tau_{0}=0$ and if $j\ge 1$, then \[\tau_{j}=\min\{i>\tau_{j-1}\mid x_{i}\in F,\ x_{i}\ne x_{\tau_{j-1}}\},\]
where $0\le j \le k$ and $k=k(\pi)$ is the last index such that the last set is nonempty. The \emph{shadow} of $\pi$ is defined as
\[\sigma(\pi)=\big(v\!=\!x_{\tau_{0}},x_{\tau_{1}},\ldots,w\!=\!x_{\tau_{k}}\big)\]
which is a path in $X_{F}$. Now let $v$ and $w$ be vertices in $F$ that are adjacent in $X_{F}$. We define
\[\Omega(v)=\sigma^{-1}(v)=\{\pi\in\Pi_{X}(v,v)\mid k(\pi)=0\}\]
and
\[\Lambda(v,w)=\sigma^{-1}\big((v,w)\big)\cap \Pi_{X}^{*}(v,w)=\{\pi\in\Pi_{X}^{*}(v,w)\mid k(\pi)=1\}.\]
For any transition matrix $Q$ of a finite directed graph we write
\[Q^* (z)=\sum_{n=0}^{\infty}(zQ)^{n}=(I-zQ)^{-1},\quad \mbox{where\quad } z\in \U.\]
The coordinates of $Q^* (z)$ are rational functions in $z$. They can easily be computed explicitly. Let
\[Q_{B}=\big(q_{B}(x,y)\big)_{x,y\in V\hat{C}}\]
denote the transition matrix of the simple random walk on $\hat{C}$ with absorbing boundary $B\subset \theta C$. This means that
\[q_{B}(x,y)=
\begin{cases}
\frac{1}{\deg_{\hat{C}} x}& \mbox{if\quad } x\sim y \mbox{\quad and\quad } x\not\in B\\
0& \mbox{otherwise.}
\end{cases}\]
(Note that some authors set $q_{B}(x,x)=1$ for vertices $x$ in the absorbing boundary to get stochastic instead of substochastic transition matrices.) Let $v$ and $w$ be two different vertices in $\theta C$. Then we define
\[d,f:\U\to \OC,\quad d(z)=\big(Q^*_{\theta C\backslash \{v\}}(z)\big)_{v,w}\cdot (\theta_X-1)\mbox{\quad and \quad} f(z)=\big(Q^*_{\theta C\backslash \{v\}}(z)\big)_{v,v}.\]
The function $d$ is independent of the choice of $v$ and $w$ because $\Aut(\hat{C})$ acts doubly transitively on $\theta C$ and $f$ is independent of $v$ by the simple transitivity of $\Aut(\hat{C})$ on $\theta C$. Now $d$ is the generating function of the probabilities that the simple random walk on $\hat{C}$ starting in some vertex $v$ in $\theta C$ hits a vertex in $\theta C\backslash \{v\}$ for the first time after $n$ steps, whereas $f$ describes the $n$-step return probabilities of the random walk to $v$ without hitting the absorbing boundary $\theta C\backslash \{v\}$. We call $d$ the \emph{transition function} and $f$ the \emph{return function}. 

\begin{lem}\label{functionaleqn}
Let $X$ be a graph which is symmetrically self-similar with respect to $F$ and $\phi$, and let $v$ and $w$ be vertices in $F$. Then
\[G(v,w|z)=G(\phi v,\phi w|d(z))\cdot f(z).\]
\end{lem}

\begin{proof}
For sets of paths $\Pi_{1}$ and $\Pi_{2}$ let $\Pi_{1}\circ\Pi_{2}$ be the set of all possible concatenations of paths in $\Pi_{1}$ with paths in $\Pi_{2}$. Then for any path
\[\pi=(v_{o}\!=\!v,v_{1},\ldots,v_{n}\!=\!w)\]
in $\Pi_{X_{F}}(v,w)$ we have
\[\sigma^{-1}(\pi)=\Lambda(v_{0},v_{1})\circ\Lambda(v_{1},v_{2})\circ\ldots\circ\Lambda(v_{n-1},v_{n})\circ\Omega(v_{n}).\]
For the weight function this implies
\[W(\sigma^{-1}(\pi)|z)=\prod_{i=0}^{n-1}W(\Lambda(v_{i},v_{i+1})|z)\cdot W(\Omega(v_{n})|z).\]
The function $d$ describes the transition from $v_{i}$ to any vertex in a boundary of any cell which is adjacent to $v_{i}$ and by Lemma \ref{cells}
\[W(\Lambda(v_{i},v_{i+1})|z)=\frac{d(z)}{\cells_{X} v_{i}\cdot(\theta_X-1)}=\frac{d(z)}{\deg_{X_{F}} v_{i}}\]
and
\[W(\Omega(v_{n})|z)=f(z).\]
Thus
\[G(v,w|z)=W(\Pi_{X}(v,w)|z)=\sum_{\pi\in\Pi_{X_{F}}(v,w)}W(\sigma^{-1}(\pi)|z)=\]
\[\sum_{\substack{
\pi\in\Pi_{X_{F}}(v,w)\\
\pi=(v_{o}=v,v_{2},\ldots,v_{n}=w)}} \prod_{i=0}^{n-1}W(\Lambda(v_{i},v_{i+1})|z)\cdot W(\Omega(v_{n})|z)=\]
\[\sum_{\substack{
\pi\in\Pi_{X_{F}}(v,w)\\
\pi=(v_{o}=v,v_{2},\ldots,v_{n}= w)}} \prod_{i=0}^{n-1}\frac{d(z)}{\deg_{X_{F}}v_{i}}\cdot f(z)=\]
\[G_{X_F}\big(v,w|d(z)\big)\cdot f(z)=G\big(\phi v,\phi w|d(z)\big)\, f(z),\]
where $G_{X_F}$ denotes the Green function on $X_F$.
\end{proof}

\begin{rem}\label{proofsteps}
The main arguments in the proof of Lemma \ref{functionaleqn} are:
\begin{enumerate}[(i)]
\item The graphs $X$ and $X_F$ are isomorphic. The set of shadows of all paths in $X$ connecting two vertices $v$ and $w$ in $F$ is the set of all paths in $X_F$ connecting $v$ and $w$.
\item A transition along an edge in $\E X_F$ from $v_i$ to $v_{i+1}$ corresponds to all possible paths in $X$ from $v_i$ to $v_{i+1}$ in the sense of the transition function $d$. These paths starting in $v_i$ may return to $v_i$ arbitrarily often and $v_{i+1}$ is the first vertex they hit in $F\backslash \{v_i\}$. This is a concatenation of generating functions in the sense of Lemma 2.2.22 in Goulden and Jackson's book \cite{goulden83combinatorial} and we substitute $z\mapsto d(z)$.
\item By the substitution in (ii) we did not consider the possibility of returning again to $w$, after hitting $w$ for the first time. The concatenation of paths that hit $w$ for the first time with their last vertex and paths that start in $w$ and return to $w$ arbitrarily often corresponds to a product of generating functions with distinct configurations in the sense of Lemma 2.2.14 in \cite{goulden83combinatorial}. Thus $G\big(\phi v,\phi w|d(z)\big)$ has to be multiplied by $f(z)$.
\end{enumerate}
\end{rem}

\section{Decomposition of the Green functions}\label{decompose}

Let $C$ be a cell of a symmetrically self-similar graph $X$. For any $x\in C$ and $y\in \overline{C}$ we define
\[h_{C}(x,y|z)=\big(Q^*_{\theta C}(z)\big)_{x,y}\]
as the $(x,y)$-coordinate of $Q^*_{\theta C}(z)$. In other words, $h_{C}(x,y|\,\cdot\,)$ is the generating function of the $n$-step transition probabilities of the simple random walk on $\hat{C}$ with absorbing boundary $\theta C$ starting in $x$ and ending in $y$. If $x\in \theta{C}$ and $y\in C$, then we define $\tilde h(x,y|z)=\frac{\deg_{\hat C}y}{\deg_{\hat C}x}h(y,x|z)$. This is the generating function for the probabilities that the simple random walk on $\hat C$, starting in $x$, hits $y$ after $n$ steps without returning to $x$ and without hitting any other vertex in $\theta C$. These functions $h$ and $\tilde h$ are called \emph{inner transition functions}. Axiom (S1) implies that there are only finitely many different inner transition functions. Note that any simple random walk is reversible and that $\deg_X a\cdot G(a,b|z)=\deg_X b\cdot G(b,a|z)$ for any vertices $a$ and $b$. The cell which contains some given vertex $x\in \V X\backslash F$ is denoted by $C(x)$.

\begin{thm}\label{algorithm}\mbox{}\par
Let $x,y\in \V X\backslash F$ and $v,w\in F$, then
\begin{enumerate}[(i)]
\item \[G(v,w|z)=G\big(\phi v,\phi w|d(z)\big)\cdot f(z),\]
\item
\[G(x,w|z)= \sum_{\substack{v\in \theta C(x)}}h_{C(x)}(x,v|z)\cdot G\big(\phi v,\phi w|d(z)\big)\cdot f(z),\]
\[G(v,y|z)= \sum_{\substack{w\in \theta C(y)}}G\big(\phi v,\phi w|d(z)\big)\cdot f(z)\cdot {\tilde h}_{C(y)}(w,y|z),\]
\item and
 \[G(x,y|z)= \delta_{C(x),C(y)}\cdot h_{C(x)}(x,y|z)\]
\[+\sum_{\substack{v\in \theta C(x)\\ w\in \theta C(y)}}h_{C(x)}(x,v|z)\cdot G\big(\phi v,\phi w|d(z)\big)\cdot f(z) \cdot {\tilde h}_{C(y)}(w,y|z),\]
where $\delta$ denotes the usual Kronecker symbol.
\end{enumerate}
\end{thm}

\begin{proof}
Let $x$ and $y$ be two vertices in $\V X\backslash F$ such that $C(x)\ne C(y)$. We decompose paths of the simple random walk on $X$ starting in $x$ and ending in $y$ into
\begin{enumerate}[(a)]
\item a path from $x$ to a vertex $v$ whose vertices are all contained in $C(x)$ except for the last vertex $v\in\theta C(x)$,
\item a path from $v$ to a vertex $w\in \theta C(y)$ and 
\item a path from $w$ to $y$ which is contained in $C(y)$ except for its first vertex $w$.
\end{enumerate}
This implies the product decomposition of (iii):
\begin{align}
G(x,y|z)&=\sum_{\substack{v\in \theta C(x)\\ w\in \theta C(y)}}h_{C(x)}(x,v|z)\cdot G(v,w|z) \cdot  {\tilde h}_{C(y)}(w,y|z)\nonumber\\
&=\sum_{\substack{v\in \theta C(x)\\ w\in \theta C(y)}}h_{C(x)}(x,v|z)\cdot G\big(\phi v,\phi w|d(z)\big)\cdot f(z) \cdot {\tilde h}_{C(y)}(w,y|z).\nonumber
\end{align}
If $x$ and $y$ lie in the same cell $C(x)=C(y)$, then paths from $x$ to $y$ need not contain vertices of $\theta C(x)$ and we get the additional additive term $h_{C(x)}(x,y|z)$.\par
Let $x$ be a vertex in $\V X\backslash F$ and let $w$ be a vertex in $F$. Then part (c) of the above path decomposition and therefore also the functions $\tilde h_{C(y)}(w,y|\,\cdot\,)$ are cancelled. This implies (ii).\par
The identity (i) corresponds to Lemma \ref{functionaleqn}.
\end{proof}

By iterating formula (i) in Theorem \ref{algorithm}, we obtain:

\begin{cor}\label{n_equation}
Let $v$ and $w$ be vertices in $F^n$. Then
\[G(v,w|z)=G(\phi^n v,\phi^n w|d^n(z))\prod_{k=0}^{n-1}f(d^k(z)).\]
For arbitrary vertices $x$ and $y$ we have
\[G(\psi^n x,\psi^n y|z)=G(x,y|d^n(z))\prod_{k=0}^{n-1}f(d^k(z)).\]
\end{cor}

\section{Green functions on $n$-cells}\label{ncells}

Let $\poles (g)$ be the set of poles and let $\zeroes (g)$ be the set of zeroes of a complex function $g$. For a cell $C$ of a symmetrically self-similar graph $X$ we define
\[\poles(\hat C)=\bigcup_{x\in C,\ y\in\overline C} \poles \big(h_{C}(x,y|\,\cdot\,)\big).\]
Axiom $(S1)$ of course implies that $\poles(\hat C)$ does not depend on the choice of the cell $C$.

\begin{lem}\label{zeroes}
The zeroes of the return function $f$ are poles of the inner transition functions $h$,
\[\zeroes (f)\subset \poles(\hat C).\]
\end{lem}

\begin{proof}
We consider the simple random walk on $\overline{C}$ with absorbing boundary $\theta C\backslash\{v\}$ for some given vertex $v$ in $\theta C$. Let $\hat{G}(x,y|\,\cdot\,)$ be the corresponding Green function on $\hat{C}$ for vertices $x$ and $y$ in $\overline{C}$ and we define
\[\hat{F}(x,y|z)=\sum_{n=0}^{\infty}q^{(n)}(x,y)z^n\]
where $q^{(n)}(x,y)$ is the $n$-step transition probability starting in $x$ for the first hit in $y$ if $n>0$, while $q^{(0)}(x,y)=0$. We have
\[\hat{F}(v,v| z)=1-\frac{1}{\hat{G}(v,v|z)}=1-\frac{1}{f(z)}.\]
A proof for this identity can be found for example in \cite[Lemma~1.13~(a)]{woess00random}. Note that any zero of $f$ is a pole of $\hat{F}(v,v|\cdot)$. For some fixed positive integer $n$, let $r^{(n)}(v,x)$ be the $n$-step transition probability from $v$ to $x$ such that the random walk meets $v$ only at the beginning. If the random walk returns to $v$ for the first time after $k$ steps, then either it returns after $k$ steps, where $2\le k\le n$, or after $n$ transitions it is in a vertex $x\in C$, without hitting $v$ during these first $n$ steps, and then it returns to $v$ (after $k\!-\!n$ steps, $1\le k\!-\!n$) in the sense of the inner transition function $h_C(x,v|\,\cdot\, )$. In other words:
\[\hat{F}(v,v|z)=\sum_{k=0}^n q^{(k)}(v,v)\cdot z^k+\sum_{x\in C}r^{(n)}(v,x)\cdot z^n\cdot h_C(x,v|z).\]
This identity holds for any positive integer $n$. In particular, any pole of $\hat{F}(v,v|\cdot)$ is a pole in $\poles(\hat C)$
which implies $\zeroes (f)\subset \poles(\hat C)$.
\end{proof}

\begin{lem}\label{green_n-cells}
Let $\hat C_n$ be the graph spanned by the closure of an $n$-cell $C_n$ and let $v$ be a vertex in the boundary $\theta C_n$. By $G_n^A(x,y|\,\cdot\,)$ we denote the Green functions of the simple random walk on $\hat C_n$ with absorbing boundary $\theta C_n\backslash\{v\}$ and by $G_n^B(x,y|\,\cdot\,)$ those with absorbing boundary $\theta C_n$. Let $w$ be a vertex in $\theta C_n\backslash\{v\}$, and let $y$ be a vertex in $F^{n-1}\cap C_n$. Then

\begin{equation}\label{eqn1}
(\theta_X-1)\cdot G_n^A(v,w|z)=d^n(z),
\end{equation}
\begin{equation}\label{eqn2}
G_n^A(v,v|z)=\prod_{k=0}^{n-1}f\big( d^k(z)\big) \mbox{\quad and}
\end{equation}
\begin{equation}\label{eqn3}
G_n^B(y,v|z)=h_{\phi^{n-1} (C_n\cap F^{n-1})}\big(\phi^{n-1} y,\phi^{n-1} v| d^{n-1}(z)\big).
\end{equation}
\end{lem}

\begin{proof}
For $n=1$ the identities \eqref{eqn1}, \eqref{eqn2} and \eqref{eqn3} are just the definitions of the functions $d$, $f$ and $h$; we recall that $F^0=\V X$. Suppose these identities are true for $n-1$. Then we use the same path decomposition and the same substitution $z\to d(z)$ as in the proof of Lemma \ref{functionaleqn}, with the difference that in \eqref{eqn1} and in \eqref{eqn3} the end point of the random walk is an absorbing vertex. Thus we do not have the possibility of returning again to the end point, which means that we do not multiply $G_{n-1}^A(v,w|z)$ and $G_{n-1}^A(v,v|z)$ with $f(z)$, see Remark \ref{proofsteps} (iii). In \eqref{eqn2} the vertex $v$ is not absorbing and we have to perform this multiplication.
\end{proof}

\section{Analytic continuation of Green functions}\label{dynamics}

Let $X$ again be a doubly symmetrically self-similar graph. The Julia set of $d$ is denoted by $\cj$.

\begin{lem}\label{fixed_point}
The point $z=0$ is a superattracting fixed point and $z=1$ is a repelling fixed point of the transition function $d$. The order of $d$ at $z=0$ is the diameter of the boundary of a cell.
\end{lem}

\begin{proof}
Let $v$ and $w$ be any two different vertices in $F$. Then $p^{(0)}(v,w)=0$. Since, by Axiom (F1), no pair of vertices in the boundary of a cell $C$ is adjacent, the order of $d$ at $z=0$ is $\diam\theta C$ which is at least 2. Thus $d(0)=0$ and $d'(0)=0$.\par
The function $d$ is probability generating which implies $d(1)=1$. Starting in a vertex $v\in F$, the derivative $d'(1)$ is the expected number of steps needed to hit any other vertex in $F$ for the first time. Hence $d'(1)>2$.
\end{proof}

We call
\[\cd=\bigcup_{n=0}^{\infty}d^{-n}(\poles (f) \cup \poles(\hat C))\]
the \emph{discrete exceptional set}.\par
Let $A_0$ be the basin of attraction of 0. Its immediate basin of attraction (the connected component of $A_0$ which contains 0) is denoted by $A^*_0$. Let $\cs$ be the set of all singularities of all Green functions.\par

\begin{thm}\label{analyticon}
There is a unique analytic continuation of any Green function to $A^*_0\backslash\cs$, which contains $A^*_0\backslash\cd$, and
\[\spec^{-1}\!P\ \subset\ \OC\backslash(A^*_0\backslash\cs)\ \subset\ \OC\backslash(A^*_0\backslash\overline\cd).\]
\end{thm}

\begin{proof}
Let $z$ be a complex number in $A^*_0$. Then there is an integer $n$ and a neighbourhood $U(z)$ of $z$ such that $U(z)\subset A^*_0$ and $d^n(U(z))$ is contained in the open unit disc $\U$. Let $x$ and $y$ be any two vertices in $\V X$. Starting with the formula for $G(x,y|\,\cdot\,)$ in Theorem \ref{algorithm} and applying it again to the Green functions $G(\phi v,\phi v|d(\,\cdot\,))$, we obtain an expression with Green functions of the form $G(\phi^2 v,\phi^2 w)|d^2(\,\cdot\,))$. Iterating this formula $n$ times in this way we get a term of the form
\begin{equation}\label{sum}
G(x,y|z)=A(z)+\sum_{\substack{\mathrm{finite}\\ \mathrm{sum}}} B(z)\cdot G(\phi^n v,\phi^n w|d^n(z))
\end{equation}
with a finite sum over Green functions of vertices $\phi^n v$ and $\phi^n w$, $v$ and $w$ are elements of $F^n$, and rational functions $A$ and $B$ which are analytic in $\OC\backslash\cd$. We obtain the analytic continuation of $G(x,y|\,\cdot\,)$ in $U(z)$ if $z$ is an element of $A^*_0\backslash\cd$. This also implies $\cs\cap A^*_0\subset\cd$. If $z$ is a point in $\cd\backslash\cs$ then again \eqref{sum} yields the analytic continuation of $G(x,y|\,\cdot\,)$ in $U(z)$. In this case $z$ is a pole in $\cd$ which is cancelled in the above sum in \eqref{sum}. Thus all Green functions are analytic in $A^*_0\backslash\cs$. Note that the choice of $n$ does only depend on $z$ and not on the vertices $x$ and $y$. A suitable choice for $n$ is such that $|d^n(z)|<1$. Uniqueness of the analytic continuation of the Green function at the point $z$ follows from the explicit form of \eqref{sum}. In particular, if $z$ is real, then we see that there is an $\varepsilon>0$ such that for any $x$ and $y$ in $\V X$ the function $G(x,y|\,\cdot\,)$ is real and analytic in $(z-\varepsilon,z+\varepsilon)$. Thus the inversion formula for the resolvent of a self adjoint operator (see Dunford and Schwarz \cite[Theorem~X.6.]{dunford63linear}) implies that $1/z$ is not in the spectrum of the transition operator $P$. This method was, for example, also used by Kesten (see \cite[Lemma~2.1]{kesten59symmetric}). The values of the Green functions at $z$ are the coordinates of the inverse operator
\[(G(x,y|z))_{x,y\in \V X}=\sum_{n=0}^\infty z^n P^n=(I-zP)^{-1}.\]
\end{proof}

\section{Non-polar singularities and the Julia set of the transition function}\label{notisolated}

Let $x$ and $y$ be vertices of a locally finite graph $Y$. We consider the simple random walk starting in $x$. Let $q^{(n)}(x,y)$ be the probability of hitting $y$ for the first time after $n$-steps. We define
\[F(x,y|z)=\sum_{n=0}^{\infty}q^{(n)}(x,y)z^n.\]
The following fact is in principle well known, though hard to be found explicitly stated in the literature.

\begin{lem}\label{lem:nonpolar}
Let the simple random walk on a locally finite, infinite graph $X$ be recurrent. Then for any pair of vertices $x$ and $y$ in $\V X$ the point $z=1$ is a non-polar singularity of the Green function $G(x,y|\,\cdot\,)$.
\end{lem}

\begin{proof}
Green functions are defined as power series around $z=0$ with positive coefficients. Pringsheim's theorem implies that the radius of convergence is a singularity. By recurrence the radius of convergence is 1. Now suppose that $z=1$ is a pole of order $k\ge 1$ of $G(x,y|\,\cdot\,)$. Then, near 1,
\[G(x,y|z)=\frac{1}{(1-z)^k}\ \tilde G(z),\]
where $\tilde G$ is analytic in a neighbourhood of 1 and $\tilde G(1)\ne 0$. We rewrite
\[G(x,y|z)=\frac{F(x,y|z)}{1-F(y,y|z)}\mbox{\qquad as\qquad}\frac{1-F(y,y|z)}{1-z}=\frac{(1-z)^{k-1} F(x,y|z)}{\tilde G(z)}.\]
Note that recurrence implies $F(x,y|1)=1$ for all pairs of vertices $x$ and $y$. Let $z\to 1^-$ along the real axis. Then we obtain $F'(y,y|1^-)<\infty$. But this is the expected return time to $y$, which has to be infinite because recurrence of the simple random walk on infinite, connected graphs is equivalent to null recurrence.
\end{proof}

We are interested in the spectrum of the Laplacian
\[\Delta=I-P\]  
on $\ell_{2}(X)$, where $I$ denotes the identity. These operators are selfadjoint and
\[\spec\Delta=\{1-\lambda\mid \lambda\in\spec P\},\quad
\spec P\subset [-1,1]\mbox{\quad and\quad } \spec \Delta\subset [0,2].\]
We recall that $\cs$ is the set of all singularities of all Green functions and that $\spec^{-1}\!P$ is the reciprocal spectrum $\{1/\lambda\mid \lambda\in \spec P\}$.\par
The next lemma follows immediately from the fact that the resolvent is analytic outside the spectrum.

\begin{lem}\label{singinspec}
For any locally finite graph we have
\[\cs\subset \spec^{-1}\!P.\]
\end{lem}

For the rest of this section let $X$ again be a symmetrically self-similar graph with bounded geometry.

\begin{lem}\label{d_real}
\[\cd\subset\OR\backslash (-1,1).\]
\end{lem}

\begin{proof}
The identities \eqref{eqn1} and \eqref{eqn3} imply that $d^{n-1}(\poles (f))$ and $d^{n-1}(\poles(\hat C))$ are poles of $G_n^A(v,w|\,\cdot\,)$ and $G_n^B(y,v|\,\cdot\,)$, respectively. Let $z_0$ be a point in $d^{n-1}(\poles (f))$. Either $z_0$ is a pole of
\[G_n^A(v,v|\,\cdot\,)=\prod_{k=0}^{n-1}f\circ d^k\]
or, if it is cancelled in this product, it is a zero of one of the factors $f\circ d^k$, for an integer $k$ with $0\le k\le n-2$. In other words, $z_0\in d^{-k}(\zeroes (f))$. Then Lemma \ref{zeroes} implies $z_0\in d^{-k}(\poles(\hat C))$ and, by the above argument, $z_0$ is a pole of a Green function $G_{k+1}^B(y,v|\,.)$. We conclude that all points in $d^{n-1}(\poles (f))$ and $d^{n-1}(\poles(\hat C))$ are poles of Green functions of type $G_n^A$ or $G_n^B$.\par
Since all our Green functions are analytic in $\U$, their singularities lie in $\OC\backslash\U$. By Lemma \ref{singinspec} the reciprocal values of these singularities are in the corresponding spectrum of a self-adjoint operator. Thus they are real and
\[d^{n-1}(\poles (f))\cup d^{n-1}(\poles(\hat C))\subset\OR\backslash (-1,1),\]
for any positive integer $n$, which implies
\[\cd\subset\OR\backslash (-1,1).\]
\end{proof}

Let us define two complex numbers $z_1$ and $z_2$ as equivalent if there are integers $m$ and $n$ such that $d^m(z_1)=d^n(z_2)$. A complex number is called $\emph{exceptional}$ if the corresponding equivalence class is finite. The Julia set of $d$ is denoted by $\cj$, the Fatou set by $\cf$.

\begin{thm}\label{nonpolar}
For the reciprocal spectrum we have
\[\spec^{-1}\!P\ =\ \cj\cup\cs\ \subset\ \cj\cup\cd.\] 
It is the set of points $z$ in $\OR\setminus (-1,1)$ such that there is a Green function which cannot be continued analytically from both half spheres of $\OC\setminus\OR$ to $z$. The super attracting fixed point $z=0$ is the only attracting fixed point of $d$, $\cf=A^*_0=A_0$ and $\cj$ is the set of accumulation points of $\cd$. Either $d$ has no exceptional points or $z=0$ is the only exceptional point of $d$. In particular, the discrete exceptional set $\cd$ does not contain exceptional points.
\end{thm}

\begin{proof}
The immediate basin of attraction $A^*_0$ of $z=0$ is the component of the Fatou set $\cf$ which contains $z=0$, see Theorm~6.3.1 in the book of Beardon \cite{beardon91iteration}. Suppose there is a point $z_0$ in $\partial A^*_0$, which is a subset of $\cj$, and a Green function $G_0$ which continues analytically from $A^*_0$ to an open neighbourhood $U(z_0)$ of $z_0$. The Julia set contains all repelling fixed points, and the backward iterates of any point in $\cj$ are dense in $\cj$, see for example \cite[Theorem~4.2.7 (ii)]{beardon91iteration} or Theorems~1.1 and 1.6 in Chapter III of the book of Carlson and Gamelin \cite{carleson93complex}. Thus there is a $\tilde z$ in $U(z_0)$ and an integer $n$ such that $d^n(\tilde z)=1$. The function $d$ is finitely branched. Corollary \ref{n_equation} implies that the Green function $G_0$ continues analytically to an open neighbourhood of $z=1$ except for a finite set of singularities. Now $z=1$ is an isolated singularity of $G_0$ which is non-polar (see Lemma \ref{lem:nonpolar}), therefore it is an essential singularity. Corollary \ref{n_equation} also implies that all $d$-backward iterates of $z=1$ are essential singularities of $G_0$. These backward iterates are dense in $\cj$ in contradiction to the assumption that $G_0$ continues analytically to an open neighbourhood $U(z_0)$ of $z_0$ in $\cj$. Since all Green functions are analytic on $\C\setminus\spec^{-1}\!P$, this implies $\partial A^*_0\subset\spec^{-1}\!P$. The reciprocal spectrum $\spec^{-1}\!P$ is real and it follows that $z=0$ is the only attracting fixed point of $d$ and $A^*_0=A_0=\cf$. This implies $\cj=\partial A^*_0$. The set $\cs\cap\cf$ is contained in $\spec^{-1}\!P$ and Theorem \ref{analyticon} implies that $\cf\setminus\cs$ is contained in $\C\setminus\spec^{-1}\!P$. We conclude,
$\spec^{-1}\!P=\cj\cup\cs$. Since $d$ is rational of degree at least 2, there are at most two exceptional points (\cite[Theorem 4.1.2]{beardon91iteration}). Exceptional points are contained in the Fatou set (\cite[Corollary~4.1.3]{beardon91iteration}). Thus $z=0$ is the only exceptional point, or there are no exceptional points. By Lemma \ref{d_real}, there are no exceptional points in $\cd$ and Theorem 4.2.7 (i) in \cite{beardon91iteration} now implies that $\cj$ is the set of accumulation points of $\cd$.
\end{proof}

Combining Theorems \ref{analyticon} and \ref{nonpolar}, we obtain:

\begin{thm}\label{thm:bounds}
\[\cj\ \subset\ \spec^{-1}\!P\ =\ \cj\cup\cs\ \subset\ \cj\cup\cd\ =\ \overline\cd\ \subset\ \OR\backslash (-1,1).\]
\end{thm}

Note that $\spec\Delta=\{1-\lambda\mid1/\lambda\in\spec^{-1} P\}$.

\section{The Julia set of the transition function}\label{julia}

From the theory of iteration of rational maps we know that the Julia set of a rational function of degree greater or equal 2 is either connected in the complex sphere $\OC$ or has uncountably many components and each point in the Julia set is a accumulation point of infinitely many distinct components of $\cj$, see for example \cite[Theorem~5.7.1]{beardon91iteration} or \cite[Lemma~4]{inninger01rational}. In our case the Julia set is real and we obtain the following theorem:

\begin{thm}
Either the Julia set of the transition function of a symmetrically self-similar graph is a real interval in $\OC$ or it is a Cantor set.
\end{thm}

\begin{exmp}
Let $o$ be an origin vertex of the two sided infinite line. The graphs $\hat C$ can be paths of length $n$ for any integer $n$. Technical calculations show that for the corresponding transition function $d_n$ we have
\[d_n(z)=z^n/\Big(\sum_{k=0}^{\lfloor \frac{n}{2}\rfloor}\binom{n}{2k} (1-z^2)^k\Big),\]
these functions are conjugated to the Chebychev polynomials $T_n$,
\[T_n(z)=\cos(n \cos^{-1}z)=1/d_n(1/z),\]
and $\cj=\OR\backslash (-1,1)$. We consider this Green function as the solution
\[\prod_{k=0}^\infty f(d_2^k(z)),\qquad \mbox{for\quad} f(z)=\frac{2}{2-z^2}\qquad\mbox{and\quad} d_2(z)=\frac{z^2}{2-z^2},\]
of the functional equation
\[G(o,o|z)=f(z)\,G(o,o|d_2 (z)),\]
where the graphs $\hat C$ spanned by the closures of the cells are paths of length 2. This solution converges as an infinite product in $\C\backslash \cj$ to
\[G(o,o|z)=\frac{1}{\sqrt{1-z^2}}, \qquad z\in\C\backslash \big((-\infty,1]\cup [1,\infty)\big).\]
The fact that the product diverges on $\cj$ corresponds to the existence of two different branches of the square root.
\end{exmp}

\begin{thm}
Let $\cj$ be the Julia set of the transition function of a symmetrically self-similar graph. Then the reciprocal spectrum $\spec^{-1}P$ coincides with the set of singularities of all Green functions $\cs$.
\end{thm}

\begin{proof}
The branches of the analytic continuations from $\U$ to the half spheres in $\OC\backslash \OR$ coincide at the real part of the Fatou set $\cf\cap\OR$ which is dense in $\cj$ because $\cj$ is supposed to be a Cantor set. None of these branches can be continued locally analytically to any point in $\cj$ since this local continuation would then have to coincide with both branches of the Green function, which is impossible since no Green function can be continued analytically to a point in the Julia set. For points in $\cd\backslash \cj$ we have already proved that they are in the spectrum if and only if they are singularities of a Green function, see Theorem \ref{analyticon} and Lemma \ref{singinspec}.
\end{proof}

For symmetrically self-similar graphs with bounded geometry we conjecture that the Julia set $\cj$ is an interval if and only if the cell-graph $\hat C$ is a finite line, otherwise it is a Cantor set. According to Lemma \ref{components}, the cell graph $\hat C$ is a finite line if and only if $X$ is a `star' consisting of finitely many one-sided infinite lines which have exactly one vertex in common.\par
Sabot showed in \cite{sabot00pure} that the integrated density of states for nested fractals with at least three essential fixed points is completely generated by the so-called Neumann-Dirichlet eigenvalues. In \cite{malozemov95pure} Malozemov and Teplyaev proved similar results for symmetric self-similar graphs where each cell has exactly two boundary points.
\\
\par\noindent
\emph{Acknowledgement.}
The author wants to thank Peter Grabner, Elmar Teufl and Wolfgang Woess for fruitful discussions. Representative for many contributions: the idea of choosing a proper set $F$ in Definition \ref{def:ssg} (self-similarity of graphs) is due to Peter Grabner, Elmar Teufl found a right idea for the proof of Lemma \ref{zeroes} and Wolfgang Woess is the author of the proof of Lemma \ref{lem:nonpolar}.

\end{document}